\newif\ifels\elsfalse
\documentclass[1p,final]{elsarticle} \elstrue

\usepackage{amssymb}
\usepackage{epsfig}
\ifels
\usepackage[justification=centering]{caption}
\fi
\usepackage{amsmath}
\usepackage{mathtools}
\usepackage{bm}
\usepackage{url}
\usepackage{multirow}
\usepackage{booktabs}

\usepackage{showlabels} 

\usepackage{verbatim}

\title{Efficient cubature rules%
}

\author{James R. Van Zandt
\ifels
\fnref{mitre}
\fi
}

\ifels
\ead{jrv@mitre.org}
\fi

\ifels
\fntext[mitre]{The MITRE Corporation, 202 Burlington Rd, Bedford, MA 01730.  The author's affiliation
  with The MITRE Corporation is provided for identification purposes
  only, and is not intended to convey or imply MITRE's concurrence
  with, or support for, the positions, opinions or viewpoints
  expressed by the author.}
\fi

\newcommand\T{\rule{0pt}{2.6ex}}      
\newcommand\TT{\rule{0pt}{3.2ex}}     
\newcommand\B{\rule[-1.2ex]{0pt}{0pt}}

\long\def\symbolfootnote[#1]#2{\begingroup%
  \def\thefootnote{\fnsymbol{footnote}}\footnote[#1]{#2}\endgroup}

\begin{document}

\begin{abstract}
  67 new cubature rules are found for three standard 
  multidimensional integrals with spherically symmetric regions and
  weight functions, using direct search with a numerical zero-finder.
  63 of the new rules have fewer integration points than
  known rules of the same degree, and 20 are within three points
  of M{\"o}ller's lower bound.
  Most have all positive coefficients and most have some symmetry,
  including some supported by one or two concentric spheres.  They
  include degree 7 formulas for integration over the sphere and
  Gaussian-weighted integrals over all space, each in 6 and 7 dimensions,
  with 127 and 183 points, respectively.
\end{abstract}

\ifels
\maketitle
\fi

\symbolfootnote[0]{\copyright2019 The MITRE Corporation. ALL RIGHTS
  RESERVED. Approved for Public Release; Distribution Unlimited. Case
  Number 17-4134.}

\section{Introduction}

We are concerned with estimating multi-dimensional 
integrals of the form
\def\aalpha{\bm{\alpha}}
\def\aa{{\bf a}}
\def\bb{{\bf b}}
\def\xx{{\bf x}}
\def\WW{{\bf W}}
\begin{equation}
  \label{eq:G}
  \int_\Omega w(\xx) f(\xx)\; d\xx,
\end{equation}
where $\xx=[x_1\;x_2\cdots x_n]^T$, for
the integration regions $\Omega$ and weighting functions $w(\xx)$ shown in
Table \ref{tab:integrals}.
Applications of \eqref{eq:G} include evaluation of quantum-mechanical matrix elements
with Gaussian wave functions in atomic physics \cite{offenhartz70orbital},
nuclear physics \cite{heyde90nuclear},
and particle physics \cite{isgur-karl78baryons}.
For applications in statistics, particularly Bayesian inference, see \cite{evans-swartz92integration}.
For applications in target tracking, see \cite{arasaratnam-haykin09cubature,jia-xin-cheng11sparse}.

\begin{table}[t]  \caption{Integrals Studied.}
  \centering
  \begin{tabular}{ccc}
    Name & Region $\Omega$ & Weight Function $w(\xx)$\\
\hline
$G_n$ & entire space $\mathbb{R}^n$ & (2$\pi)^{-n/2} e^{-\xx^T\xx/2}$ \T\\
$E_n^{r^2}$ & entire space $\mathbb{R}^n$ &  $e^{-\xx^T\xx}$ \\
$E_n^r$ & entire space $\mathbb{R}^n$ &  $e^{-\sqrt{\xx^T\xx}}$ \\
$S_n$ & unit $n$-sphere $\xx^T\xx \le 1 $ &  1 \\
  \end{tabular}
  \label{tab:integrals}
\end{table}

We approximate these integrals using {\em cubature formulas} or
{\em integration rules} of the form
\begin{equation}\label{eq:Wx}
  \sum_{i=1}^N W_i f(\xx_i),
\end{equation}
where the {\em weights} $W_i$ and {\em nodes} or {\em points} $\xx_i$ are independent of
the function $f$.  

The first two integrals in the table are of course closely related.  Given an approximation of
$G_n$ of the form \eqref{eq:Wx}, we can construct an equivalent
approximation $E_n^{r^2}\approx \sum_{i=1}^N B_if(\bb_i)$ where
$\bb_i=\xx_i/\sqrt{2}$ and
$B_i=\pi^{n/2}W_i$.
In this paper we address $E_n^{r^2}$, following the numerical analysis convention.
However, in the supplemental material we quote the parameters for
the corresponding $G_n$ formulas for the convenience of
researchers using the other convention.

If an integration rule is exact for all polynomials up to and including degree
$d$, but not for some polynomial of degree $d+1$, then we say the
rule has {\em algebraic degree of exactness} (or simply {\em degree}) $d$.

One can construct cubature formulas exact for a space of polynomials
by solving the large system of polynomial equations associated with
it.  In describing this method, Cools stated that ``it is essential to
restrict the search to cubature formulas with a certain
structure'' \cite{cools02advances}.
For example, in \cite[``CUT4'' formulas]{adurthi-singla-singh12cut}, points were assumed to take the form
\begin{equation}
  \label{eq:cut4b}
  \begin{array}{rrrrr@{}lcc}
(0,      &0,     &\cdots&0,     &0&)     &W_0&1\\
(\pm\eta,&0,     &\cdots&0,     &0&)_S   &W_1&2n\\
(\pm\nu, &\pm\nu,&\cdots&\pm\nu,&\pm\nu&)&W_2&2^n
\end{array}
\end{equation}
where the notation $(\cdots)_S$ indicates that all points
obtained from these by permutation of coordinates are included, and
are assigned the same weight.
The last column gives the number of points.
This point set is {\em fully symmetric}; i.e., closed under all
coordinate permutations and sign changes.  However, relaxing this symmetry requirement can allow us to
find formulas with fewer points \cite{lyness65evaluations,wissmann-becker86cubature}.  For
example (as shown in Figure \ref{fig:plot3e2}), in two dimensions,
there is a formula of degree five with points at the vertices of a regular hexagon \cite[formula
V]{stroud-secrest63spherically}, which is closed under sign
permutations but not coordinate permutations. There is also a formula of degree 4 with
points at the vertices of a regular pentagon, which is closed under
sign changes in $x_1$ but not $x_2$; i.e., bilateral symmetry.

The objective of this work was to test whether the continuing
improvements in computer processing have made the 
``brute force''
approach---using a numerical zero-finder to solve the moment
constraint equations directly, making few or no assumptions about the
form of the points---feasible for interesting problems.
We find that, for rules with up to approximately 3000 free
parameters\footnote{Each of the $N$ points has $n$ coordinates and a
  weight, so a rule has $(n+1)N$ free parameters.}, it is no longer necessary to assume at the
outset that the points
have a particular structure.
Relieving those assumptions has made it possible to discover rules with fewer points than known rules of the
same degree, including twenty rules that come within three points of
the lower bound found by M{\"o}ller \cite{moller79nodes,lu-darmofal04gaussian}.

Section \ref{searching} describes our search method.  
Even a symmetric rule found by this method will have a random orientation,
making it inconvenient to present and use.  Section \ref{rotations}
describes how the description of a rule can be simplified by orienting
it to take best advantage of any symmetries.  
Many of the rules found have little or no symmetry.  Section \ref{pruning}
describes several procedures that can improve such a rule.  
Sections \ref{deg2}
through \ref{deg9} present the new rules that have significant symmetry.  Section
\ref{summary} lists and discusses all the new rules.  The supplemental material includes tables in double and
quad precision of all the new rules, including those with little or no
symmetry.

\section{Searching}\label{searching}

An approximation is exact for all polynomials with degree $\le d$  if
it is exact for all monomials
\begin{equation}
  \label{eq:mono}
f(\xx)=  x_1^{\alpha_1} x_2^{\alpha_2} \cdots x_n^{\alpha_n},\quad 0 \le \alpha_1
  +\cdots+ \alpha_n \le d ,
\end{equation}
where the $\alpha_i$ are all nonnegative integers.  If any of the
$\alpha_i$ are odd, then the monomial integral is zero for any of our
problems.  Let $\aalpha=[\alpha_1 \; \alpha_2 \cdots \alpha_n]^T$, where all $\alpha_i$ are even, and $\beta_i=(\alpha_i+1)/2$.  Then the monomial
integral is
\cite{folland01sphere,stroud-secrest63spherically,stroud71approximate}
\begin{equation}
  \label{eq:In}
  \begin{array}{l@{\;}l}
    I(\aalpha)&\equiv\int_\Omega w(\xx) x_1^{\alpha_1} x_2^{\alpha_2} \cdots
    x_n^{\alpha_n} \;dx_1\;dx_2\;\cdots\;dx_n\\
    \\
              &=\left\lbrace
  \begin{array}{rl}
\Gamma(\beta_1)\Gamma(\beta_2)\cdots\Gamma(\beta_n)    &\hbox{for $E_n^{r^2}$}\\
{\displaystyle 2 (\alpha_1+\cdots+\alpha_n+n-1)!\over
    \displaystyle\Gamma(\beta_1+\cdots+\beta_n)}
    \Gamma(\beta_1)\Gamma(\beta_2)\cdots\Gamma(\beta_n) 
                                                       &\hbox{for $E_n^r$}\\
{\displaystyle 1 \over \displaystyle
\Gamma(\beta_1+\cdots+\beta_n+1)}
\Gamma(\beta_1)\Gamma(\beta_2)\cdots\Gamma(\beta_n)    &\hbox{for $S_n$}\quad.\\
  \end{array}\right.
  \end{array} 
\end{equation}

For example, a rule of degree 4 for $E_2^{r^2}$ must
satisfy these 15 constraints, where $x_{ij}$ is the $j$th coordinate
of the $i$th point and each sum is over $i=1\ldots N$:
  \begin{equation}\setlength\arraycolsep{2pt}
    \begin{array}{lllll}
      \sum W_i=\pi      & \sum W_ix_{i1}=0  &  \sum W_ix_{i1}^2=\pi/2&  \sum W_ix_{i1}^3=0&  \sum W_ix_{i1}^4=3\pi/4\\[1ex]
      \sum W_ix_{i2}=0  & \sum W_ix_{i1}x_{i2}=0  &  \sum W_ix_{i1}^2x_{i2}=0&  \sum W_ix_{i1}^3x_{i2}=0  \\[1ex]
      \sum W_ix_{i2}^2=\pi/2& \sum W_ix_{i1}x_{i2}^2=0  &  \sum W_ix_{i1}^2x_{i2}^2=\pi/4  \\[1ex]
      \sum W_ix_{i2}^3=0& \sum W_ix_{i1}x_{i2}^3=0   \\[1ex]
      \sum W_ix_{i2}^4=3\pi/4  
    \end{array}
  \end{equation}
The number
of constraints increases rapidly with $n$ and $d$.  The rules in
Section \ref{r7_183_7} of dimension and degree 7 satisfy 3432
constraints.


Stroud \cite{stroud60quadrature} showed that if there is an $N$ point formula in $n$ dimensions
of degree $d$, then
\begin{equation}
  N \ge 
    {n + \lfloor d/2 \rfloor \choose \lfloor d/2 \rfloor} ,\label{eq:stroud}
\end{equation}
M{\"o}ller improved this bound for odd
degrees \cite{moller79nodes,lu-darmofal04gaussian}.
Let $d=2s-1$, then
\begin{equation}\label{eq:mlb}
  N \ge N_{MLB} \equiv\left\lbrace
    \begin{array}{ll}
      {n+s-1 \choose n} + \sum_{k=1}^{n-1} 2^{k-n} {k+s-1 \choose k} & s \hbox{ even} \\[1ex]
      {n+s-1 \choose n} + \sum_{k=1}^{n-1} (1-2^{k-n}){k+s-2 \choose k} & s \hbox{ odd}\quad.
      \end{array}
    \right.
\end{equation}

However, a formula satisfying the bound exactly may not exist.  We searched
for the rule of a given degree with the fewest points, using a binary
search between M{\"o}ller's lower bound and the number of points in a
known formula of the given degree or of the next higher degree.

We initialized each search with normally distributed points, assigning
initial weights of
\begin{equation}
  \WW=e^{-\sqrt{\xx^T\xx}},
\end{equation}
but then normalizing them so they sum to V:
\begin{equation}
  \label{eq:V}
      V\equiv\int_\Omega w(\xx)  \;dx_1\;dx_2\;\cdots\;dx_n
=\left\lbrace
  \begin{array}{rl}
\pi^{n/2}    &\hbox{for $E_n^{r^2}$}\\
{\displaystyle2(n-1)!\pi^{n/2} \over \displaystyle\Gamma(n/2)} &\hbox{for $E_n^r$}\\[1.5ex]
{\displaystyle 2  \pi^{n/2}\over \displaystyle n\Gamma(n/2)}&\hbox{for
    $S_n$}\quad,
  \end{array}\right.
\end{equation}
so the zeroth order constraint was satisfied exactly.  The points were then linearly
scaled so the second order constraints were also satisfied.

In most cases, the number of equations and unknowns were unequal (with
almost all problems becoming over-determined before $N$ reached M{\"o}ller's lower bound), so many of the
methods developed for solving nonlinear equations could not be applied.
We used {\tt fsolve} from the MATLAB
Optimization Toolbox \cite{coleman-li96minimizing,gill-murray-wright81optimization,matlab08optimization}, or {\tt UDL}\footnote{
We revised {\tt UDL} by adding a stopping criterion:
If, after any seven consecutive steps, the norm of the residual has
decreased by less than seven percent, then the search is deemed a failure.} by
Simonis \cite{simonis06inexact}.
After a failed search, we restarted with a new set of points.
After a success, we tried dropping
low-weight points, combining points with very near neighbors, and
simply restarting with fewer points.

In a few cases, an elegant formula was found easily---with the weights
on any extra points reduced to zero.  However, ordinarily with extra
points, and often even with no extra points, there are enough extra
degrees of freedom that any symmetry is lost, and successive searches
would find substantially different formulas.  In those cases, 
additional constraints were added --- moment
constraints of the next higher degree, starting with
\begin{equation}
  \label{eq:extras}
  \sum_{i=1}^N W_i x_{ij}^{d+1} = I([d+1 \; 0 \ldots 0 ]^T) ,
\end{equation}
for $1 \le i \le k \le n$.

Parameters of all searches were logged, along with the results of all
successful searches; so the results of a lucky random starting point would not be
lost.  After finding a rule for one of the integrals, we also
searched for similar rules for each of the other integrals, starting with the
same point layout and relative weights, but normalizing the weights and
scaling those points so its zeroth and second degree constraints were
satisfied exactly.

\section{Presentation}\label{presentation}
\subsection{Rotations}\label{rotations}
One a rule is found, the first step is to sort the points by radius.
If the rule is symmetric, as evidenced by several points at the
same radius and with equal weights, it is desirable to determine its structure and
if possible to express it in a simple form.  Any orthogonal
transformation of a set of points yields an equivalent set of points.
Note that any orthogonal transformation can be expressed in terms of a
skew-symmetric matrix via the Cayley transform
\cite{cayley43rotation}.
Thus, in $n$ dimensions we have $n(n-1)/2$ free parameters we can use
to orient a rule.  Often we concentrate on the sphere
supporting the fewest points, and we want to rotate to put one of those
points on the first coordinate axis.  Choosing  points in that
shell by increasing angular distance from that first point, we rotate to put a second point in the plane
defined by the first two axes, a third point in the subspace defined
by the first three axes, etc.
We call this ``aligning the axes'' to the chosen points.
Our process is as follows:

Assume we have chosen $n$ points.
Reorder the rows of the point matrix so those rows appear in
order at the top, forming a $n\times n$ submatrix we will call $A$.
The remainder of the rows form a submatrix we will call
$B$.  Use the QR decomposition to factor the transpose of $A$, so
that
\begin{equation}
  \label{eq:factor}
  A^T=RU,
\end{equation}
where $R$ is orthogonal and $U$ is upper triangular.
Taking the transpose of both sides, we have
\begin{equation}
  \label{eq:xpose}
  A=(RU)^T=U^T R^T,
\end{equation}
and right multiplying by $R$ we have
\begin{equation}
  \label{eq:tri}
  AR=U^T R^T R = U^T.
\end{equation}
Thus, right multiplying our original point matrix by $R$ gives us
\begin{equation}
  \left[ \begin{array}{c}A\\B\end{array}\right] R =   \left[ \begin{array}{c}U^T\\C\end{array}\right].
\end{equation}
$U^T$ is lower triangular.  In its first row, only the first element
is nonzero, so it represents a point along the first coordinate.  In
the second row, only the first two elements are nonzero, etc.  
This satisfies the requirements set out above.

If a rule in $n$ dimensions has $n+1$ points at the same radius (such
as the 6 inner points in the 5 dimensional rules of Section \ref{rule22}), they
typically appear at the vertices of a regular $n$ simplex.  In that case, a simple
description can be found by rotating one point to be equidistant
from all coordinate axes, with each of the other $n$ points in the
plane defined by that first point and one of
the coordinate axes.  

\subsection{Closed form expressions}\label{closed}
If a rule has enough
symmetry, we attempt to express its points and weights in closed
form.  In some cases they are integers, simple fractions, or square
roots of simple fractions, which can be identified by converting them
to a simple continued fraction and looking for a repeating pattern \cite{booher11continued}.
To guard against mathematical coincidences, our next step is to use Maxima \cite{racine-li08maxima,toth08maxima}
to confirm that the resulting rule satisfies the moment constraint equations
exactly, or (if Maxima could not simplify some expressions) with absolute
error less than $10^{-55}$.

We also attempted to identify the points as vertices of some known
polytope.  This allows the points to be described economically and
calculated directly.

\subsection{Economizing}\label{pruning}
A rule lacking any symmetry was assumed to have extra degrees of
freedom, even if it had the minumum number of points.
We tried improving its symmetry by
projecting the innermost or outermost few points to (or toward) the same radius, 
giving them all the same weight, and using that revised configuration
to start a new search.
In a few cases, this enabled us to eliminate negative weights or drop
some points.

We similarly tried to impose bilateral symmetry.  
We reoriented the rule so 
eigenvectors of the covariance of the unweighted points 
were aligned with the coordinate axes.
We then tested whether the rule was close to bilaterally symmetric
with respect to any of the axes.  If so, we searched for a similar
symmetric rule.

Initially, we simply sorted 
the points along that coordinate value, pairing the first with the
last, the second with the next to last, etc.  For each pair, we moved
the points to be exactly symmetric (i.e., moving each to the mean of
its own position and that of the reflection of the associated point) and gave each the mean of the two
original weights.  Points near the symmetry plane were moved to that plane.
This simple method failed when the initial positions were far enough
away from symmetry (e.g., if the eigenvectors did not match the
symmetry plane well enough) to change the ordering of the points along
the chosen axis.

Eventually we switched to treating the association of the original and
reflected points as a linear assignment problem, solving it with the Jonker-Volgenant-Castanon (JVC) assignment algorithm 
\cite{jonker-volgenant87assignment,malkoff97assignment},
This had the additional benefit of eliminating the need for a
threshold for being ``near'' the symmetry plane.  If a point were
assigned to its own reflection, then its adjusted position would automatically be
on the symmetry plane.

\section{Degree 2 rules}\label{deg2}

The classical second degree rules have $n+1$ points at the vertices of
an $n$ simplex.  They are usually presented in the form
\cite{stroud71approximate,julier03simplex,chun-ping-xianci07novel}

\begin{equation}
  \chi=\sqrt{n+1}\left[\begin{array}{ccccc} 
 \sqrt{1\over1\cdot2}&  \sqrt{1\over2\cdot3}& \sqrt{1\over3\cdot4}& \cdots & \sqrt{1\over n(n+1)}\\
-\sqrt{1\over1\cdot2}&  \sqrt{1\over2\cdot3}& \sqrt{1\over3\cdot4}& \cdots & \sqrt{1\over n(n+1)}\\
             0 & -\sqrt{2\over3}&\sqrt{1\over3\cdot4}& \cdots & \sqrt{1\over n(n+1)}\\
             0 &              0 &-\sqrt{3\over4} & \cdots & \sqrt{1\over n(n+1)}\\
        \vdots &          \vdots&         \vdots & \ddots & \sqrt{1\over n(n+1)}\\              
             0 &               0&              0 & \cdots & -\sqrt{n\over n+1}
  \label{eq:chun}
\end{array}\right],
\end{equation}
where each row represents a point, and the weight on each point is $1/(n+1)$.
Fan and You noticed that in three dimensions the points can be
expressed in the much simpler form \cite{fan-you09sigma}
\begin{equation}
  \label{eq:fan}
  \chi=\left[\begin{array}{rrrr}
               1& 1& 1 \\
               1&-1&-1 \\
               -1& 1&-1 \\
               -1&-1& 1 
    \end{array}\right].
\end{equation}
This can be generalized to other dimensions, yielding the set of points
\begin{equation}
  \label{eq:higher}
  \chi=\left[\begin{array}{cccccc}
               1     &     1&     1& \cdots& 1      \\     
               a     &     b&     b& \cdots& b      \\     
               b     &     a&     b& \cdots& b      \\     
               b     &     b&     a& \cdots& b      \\     
               \vdots&\vdots&\vdots& \ddots& \vdots \\
               b     &     b&     b& \cdots&a        
    \end{array}\right],
\end{equation}
with the two solutions
\begin{equation}
a = {-1+{\left(n-1\right)\, \sqrt{n+1}}\over{n}}, \qquad
b = {-1-{\sqrt{n+1}}\over{n}} ,
\end{equation}
or
\begin{equation}
a = {-1-{\left(n-1\right)\,\sqrt{n+1}}\over{n}} ,\qquad
b = {-1+{\sqrt{n+1}}\over{n}} .
\end{equation}
When a rule has $n+1$ points at the vertices of a regular simplex, it
can be rotated into one of these orientations.

%

\section{Degree 4 rules}

\subsection{Degree 4, dimension 3, 10 point rules}
\begin{figure}
  \centering
  \includegraphics[width=3in,angle=-90]{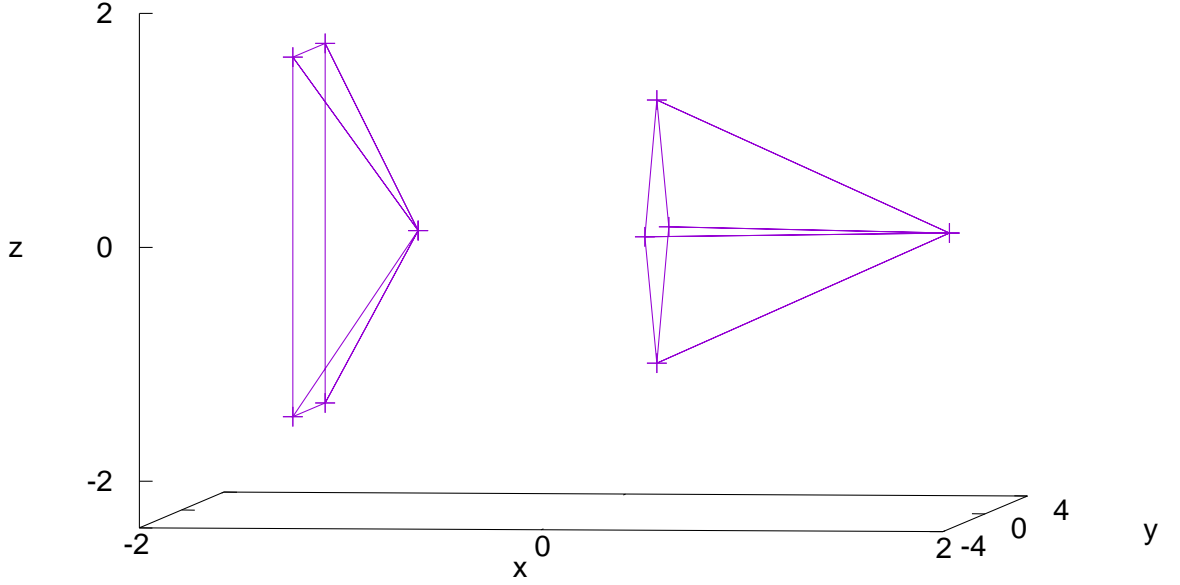}%
  \caption{Configuration of the 10 points for the rule of degree 4 for $E_3^{r^2}$.}
  \label{fig:pyramids}
\end{figure}

The points in these new formulas for $E_3^{r^2}$ and $S_3$ are closed with respect to sign changes
along two of the three coordinates.  The points form two
pyramids, with one offset and rotated from the other, as shown in
Figure \ref{fig:pyramids}.
The configuration is shown in Table \ref{tab:3_10_4}.
The Maxima program verifying the correctness of these rules is shown
in Figure \ref{fig:3_10_4_mac}.
The supplemental material includes similar programs for the other new
rules.

Becker found an 11 point cubature formula of degree 4 for $S_3$
\cite{becker87finite}, but we are not aware of any previous formulas of
degree 4 for $E_3^{r^2}$.
\begin{table}
 \centering
\caption{ 10 Point Rules of Degree 4 for $E_3^{r^2}$ and $S_3$.}\label{tab:3_10_4} 
\begin{tabular}{r@{}cr@{}cr@{}ccc}
 $ x_{1}$& &$x_{2} $& &$x_{3} $&    &Weight&  \# Points\\
\hline
 $ g    $& &  0     & &  0     &   &$ W_3$&  1 \\
 $ a    $& &$(\pm c$& &  0     &$)_S$&$ W_2$&4 \\
 $-b    $& &  0     & &  0     &   &$ W_1$&  1 \\
 $-e    $& &$\pm f $& &$\pm f $&   &$ W_4$&  4 \\
\end{tabular}
\\

\begin{tabular}{cll}\\
        & \hspace{3em} $   E_3^{r^2}\T $   & \hspace{3em} $   S_3   $\\ 
\hline
$  a   $&$(\sqrt{3}-1)/2    \T           $ &$   (2\sqrt{3}-1)/\sqrt{77}           $ \\
$  b   $&$(\sqrt{7}-1)/2                 $ &$   (2\sqrt{203}-\sqrt{77})/35        $\\
$  c   $&$ \sqrt{3-\sqrt{3}}             $ &$   \sqrt{(48-8\sqrt{3})/77}          $ \\
$  e   $&$(\sqrt{3}+1)/2                 $ &$   (2\sqrt{3}+1)/\sqrt{77}           $ \\
$  f   $&$ \sqrt{(\sqrt{3}+3)/2}         $ &$   \sqrt{(24+4\sqrt{3})/77}          $ \\
$  g   $&$(\sqrt{7}+1)/2                 $ &$   (2\sqrt{203}+\sqrt{77})/35        $\\
$  W_1 $&$\pi^{3/2}\,(2\,\sqrt{7}+7)/42  $ &$  \pi(841+32\sqrt{11}\sqrt{29})/5220 $ \\
$  W_2 $&$\pi^{3/2}\,(\sqrt{3}+2)/24     $ &$  7\pi(13+4\sqrt{3})/720             $ \\
$  W_3 $&$\pi^{3/2}\,(7-2\,\sqrt{7})/42  $ &$  \pi(841-32\sqrt{11}\sqrt{29})/5220 $ \\
$  W_4 $&$\pi^{3/2}\,(2-\sqrt{3})/24     $ &$  7\pi(13-4\sqrt{3})/720             $ \\
\end{tabular}
\end{table}

\begin{figure}
  {\footnotesize
\begin{verbatim}
rule:"3_10_4"; n:3; N:10; deg:4; /* 3 dimensions, 10 points, degree 4 */    
/*            G_n                        E_n^{r^2}                S_n */
a_vals:[(sqrt(6)-sqrt(2))/2 ,(sqrt(6)-sqrt(2))/sqrt(8) ,0,(2*sqrt(3)-1)/sqrt(77)         ];
b_vals:[(sqrt(14)-sqrt(2))/2,(sqrt(14)-sqrt(2))/sqrt(8),0,(2*sqrt(203)-sqrt(77))/35      ];
c_vals:[sqrt(6 - 2*sqrt(3)) ,sqrt((6 - 2*sqrt(3))/2)   ,0,sqrt((48-8*sqrt(3))/77)        ];
e_vals:[(sqrt(6)+sqrt(2))/2 ,(1+sqrt(3))/2             ,0,(2*sqrt(3)+1)/sqrt(77)         ];
f_vals:[sqrt(sqrt(3)+3)     ,sqrt((sqrt(3)+3)/2)       ,0,sqrt((24+4*sqrt(3))/77)        ];
g_vals:[(1+sqrt(7))/sqrt(2) ,(1+sqrt(7))/2             ,0,(2*sqrt(203)+sqrt(77))/35      ];
w1vals:[(7+2*sqrt(7))/42,%pi^(3/2)*(2*sqrt(7)+7)/42,0,%pi*(841+32*sqrt(11)*sqrt(29))/5220];
w2vals:[(sqrt(3)+2)/24  ,%pi^(3/2)*(sqrt(3)+2)/24  ,0,%pi*(13+4*sqrt(3))*7/720           ];
w3vals:[(7-2*sqrt(7))/42,%pi^(3/2)*(7-2*sqrt(7))/42,0,%pi*(841-32*sqrt(11)*sqrt(29))/5220];
w4vals:[(2-sqrt(3))/24  ,%pi^(3/2)*(2-sqrt(3))/24  ,0,%pi*(13-4*sqrt(3))*7/720           ];
plabel:["       Gn","E_n^{r^2}","      E_n","       Sn"];
for problem in [1,2,4] do block(
a: a_vals[problem],
b: b_vals[problem],
c: c_vals[problem],
e: e_vals[problem],
f: f_vals[problem],
g: g_vals[problem],
w1:w1vals[problem],
w2:w2vals[problem],
w3:w3vals[problem],
w4:w4vals[problem],
x:matrix([ g, 0, 0, w3],
         [ a, c, 0, w2],
         [ a,-c, 0, w2],
         [ a, 0, c, w2],
         [ a, 0,-c, w2],
         [-b, 0, 0, w1],
         [-e, f, f, w4],
         [-e, f,-f, w4],
         [-e,-f, f, w4],
         [-e,-f,-f, w4]),
ex(a,p) := if p=0 then 1 else a^p,
alpha : makelist(0,i,n),
pass : true,
for i1 : 0 thru deg do 
  (alpha[1] : i1,
  for i2 : 0 thru deg-i1 do 
    (alpha[2] : i2,
    for i3 : 0 thru deg-i1-i2 do 
      (alpha[3] : i3,
        if sum(alpha[i],i,1,n)<=deg then 
        (beta : (alpha+makelist(1,i,n))/2,
          odd : false, for i : 1 thru n do odd : odd or oddp(alpha[i]),
  m :  if odd then 0 else 
       if problem=1 then 2^sum(beta[ii],ii,1,n)*product(gamma(beta[ii]),ii,1,n)
         /(2*%pi)^(n/2) else
       if problem=2 then product(gamma(beta[ii]),ii,1,n) else 
       if problem=3 then 2*gamma(sum(alpha[ii],ii,1,n)+n)
         *product(gamma(beta[ii]),ii,1,n)/gamma(sum(beta[ii],ii,1,n)) else 
       if problem=4 then product(gamma(beta[ii]),ii,1,n)/sum(beta[ii],ii,1,n)
         /gamma(sum(beta[ii],ii,1, else 
       if problem=5 then 2*product(gamma(beta[ii]),ii,1,n)
          /gamma(sum(beta[ii],ii,1,n)),    
          value : radcan(sum(product(ex(x[i,j],alpha[j]),j,1,n)*x[i,n+1],i,1,N)),
          subtest : value=m, /* test one constraint */
          pass : pass and subtest)))),
print(plabel[problem],":",rule,"...", if pass then "pass" else "FAIL" ))$
\end{verbatim}}
  \caption{The Maxima program ({\tt m3\_10\_4.mac}) verifying the correctness of the rules of degree 4 in 3 dimensions.}
  \label{fig:3_10_4_mac}
\end{figure}

\subsection{Degree 4, dimension 3, 11 point rule}

We were unable to find a 10 point rule for $E_3^r$, but we did find the 11
point rule shown in Table \ref{tab:3_11_4}.
This rule has at least one remaining degree of freedom, as the $x_3$ coordinate
in the fourth line of the table need not be zero.  An example with a
nonzero value appears in the supplemental material.
\begin{table}
\caption{11 Point Rule of Degree 4 for $E_3^r$.}\label{tab:3_11_4}
\begin{center}\begin{tabular}{rrrr}
$ x_{1}$& $ x_{2}$& $ x_{3}$& Weight \\
\hline
$\pm$5.123512671436 &  4.925613098468  &  0.000000000000  & 0.379658096396  \\
$\pm$4.102816292737 & -1.218122471265  &  1.544992698170  & 1.815112382679  \\
$\pm$3.636092685910 & -1.218122471265  & -4.843920857272  & 0.737101279022  \\
     0.000000000000 & -1.836923221948  &  0.000000000000  & 8.813498359176  \\
     0.000000000000 & -12.639707409137 & -3.423767380484 &  0.036648025338 \\
     0.000000000000 &  1.948389609086  & -1.422580596634  & 7.054048788228  \\
     0.000000000000 &  1.703608086180  &  3.398957047139  & 3.331366718822  \\
     0.000000000000 & -8.635010968135  & 11.051160549267  & 0.033435820963  \\

\end{tabular}
\end{center}
\end{table}

\subsection{Degree 4, dimension 4, 16 point rules} 

We found two sets of 16 point rules.  Those in the first set have a central point, one shell
of 10 points, and another shell of 5 points, as shown in Table \ref{tab:4_16_4a}.
The formula for $S_4$ has five points outside the region, and one with a negative weight.
\begin{table}
\caption{16 Point Rules of Degree 4 in 4 Dimensions (Group 1).}\label{tab:4_16_4a}
\begin{center}\begin{tabular}{r@{}cr@{}cr@{}cr@{}cccccc}
$ x_{1}$& &$ x_{2}$& &$ x_{3}$& &$ x_{4}$& &Weight & Radius &  \# Points\\
\hline
   0  & &   0  & &   0  & &   0  & &$ W_0$&0& 1 \\
$  (c $& &$ c $& &$ -b $& &$ -b $&$)_S$&$ W_1$&$r_1$& 6 \\
$ (-e $& &$ -a $& &$ -a $& &$ -a $&$)_S$ &$ W_1$&$r_1$& 4 \\
$  f $& &$  f $& &$  f $& &$  f $& &$ W_2$&$r_2$&     1 \\
$  (g $& &$ -e $& &$ -e $& &$ -e $&$)_S$&$ W_2$&$r_2$& 4 \\
\end{tabular}

\begin{tabular}{clll}\\ 
       &\hspace{3em} $   E_4^{r^2}\T $&\hspace{3em} $   E_4^r   $&   \hspace{3em} $   S_4    $\\
\hline
$ a   $&$(3\sqrt{3}-\sqrt{15})/12\T$&$(3\sqrt{42}-\sqrt{210})/12$&$ (3\sqrt{3}-\sqrt{15})/24 $\\
$ b   $&$ ( \sqrt{15}-\sqrt{3})/6  $&$(\sqrt{210}-\sqrt{42})/6   $&$ ( \sqrt{15}-\sqrt{3})/12 $\\
$ c   $&$ ( \sqrt{15}+\sqrt{3})/6  $&$(\sqrt{210}+\sqrt{42})/6   $&$ ( \sqrt{15}+\sqrt{3})/12 $\\
$ e   $&$ ( \sqrt{15}+\sqrt{3})/4  $&$(\sqrt{210}+\sqrt{42})/4   $&$ ( \sqrt{15}+\sqrt{3})/8  $\\
$ f   $&$ \sqrt{3}                 $&$\sqrt{42}                 $&$ \sqrt{3}/2               $\\
$ g   $&$ (3\sqrt{15}-\sqrt{3})/4  $&$(3\sqrt{210}-\sqrt{42})/4 $&$ (3\sqrt{15}-\sqrt{3})/8  $\\
$ W_0 $&$ \pi^2/12                  $&$29\pi^2/7               $&$ -\pi^2/9                $\\
$ W_1 $&$ 9\pi^2/100                $&$27\pi^2/35              $&$ 3\pi^2/50               $\\
$ W_2 $&$ \pi^2/300                 $&$\pi^2/35                 $&$ \pi^2/450               $\\
$ r_1 $&$ \sqrt{2}                  $&$\sqrt{28}                $&$ \sqrt{1/2}              $\\
$ r_2 $&$ \sqrt{12}                 $&$\sqrt{168}               $&$ \sqrt{3}                $\\
\end{tabular}
\end{center}
\end{table}

Each rule in the second set has a central point, a shell of 6 points,
and another shell of 9 points, as shown in Table \ref{tab:4_16_4b}.
The second formula for $S_4$ has zero weight on the central point,
making it a 15 point formula, with all positive weights and nine
points on the boundary.

Both of these formulas for $S_4$ are distinct from the 16 point formula
found by Mysovskih \cite{cools99cubature,mysovskih92cubature} which
has positive weights and seven points on the boundary.  We are
not aware of previous degree 4 formulas for the other integrals.
\begin{table}
\caption{16 or 15 Point Rules of Degree 4 in 4 Dimensions (Group 2).}\label{tab:4_16_4b}
\begin{center}\begin{tabular}{r@{}cr@{}cr@{}cr@{}ccccc}
$ x_{1}$  & &$ x_{2}$ & &$ x_{3}$& &$ x_{4}$& &Weight& Radius& \# Points\\
\hline
   0     & &   0     & &   0  & &   0  & &$ W_0$ &  0    & 1 \\
   0     & &   0     & &   0  & &$ -c $& &$ W_1$ & $r_1$ & 1  \\
   0     & &   0     & &$  c $& &   0  & &$ W_1$ & $r_1$ & 1  \\
$\pm b $ & &   0     & &$ -a $& &   0  & &$ W_1$ & $r_1$ & 2  \\
   0     & &$\pm b $ & &   0  & &$  a $& &$ W_1$ & $r_1$ & 2  \\
$\pm b $ & &$\pm b $ & &$  a $& &$ -a $& &$ W_2$ & $r_2$ & 4  \\
   0     & &$\pm b $ & &$ -c $& &$ -a $& &$ W_2$ & $r_2$ & 2  \\
$\pm b $ & &   0     & &$  a $& &$  c $& &$ W_2$ & $r_2$ & 2  \\
   0     & &   0     & &$ -c $& &$  c $& &$ W_2$ & $r_2$ & 1  \\
\end{tabular}

\begin{tabular}{lccc}\\
                &$   E_4^{r^2}$\T  &  $  E_4^r $  &  $   S_4 $ \\
\hline
$          a   $&$\sqrt{1/2}$\T &$\sqrt{7}$    &$\sqrt{1/8}$\\ 
$          b   $&$\sqrt{3/2}$   &$\sqrt{21}$   &$\sqrt{3/8}$\\ 
$          c   $&$\sqrt{2}$     &$\sqrt{28}$   &$\sqrt{1/2}$\\ 
$          W_0 $&$\pi^2/4$      &$39\pi^2/7$   &0           \\ 
$          W_1 $&$\pi^2/12$     &$5\pi^2/7$    &$\pi^2/18$  \\ 
$          W_2 $&$\pi^2/36$     &$5\pi^2/21$   &$\pi^2/54$  \\ 
$        c=r_1 $& $\sqrt{2}$    &$\sqrt{28}$   &$\sqrt{1/2}$\\ 
$          r_2 $&     2         & $\sqrt{56}$  &1           \\   
\end{tabular}
\end{center}
\end{table}

\subsection{Degree 4, dimension 5, 22 point rules}\label{rule22}
Each of these new formulas has a central
point and two shells.
We can describe the points using five generators,
as shown in Table \ref{tab:5_22_4_p}.
The six points
with weight $W_0$ are at the vertices of a regular 5-simplex.
The 15 points with weight $W_1$ are the vertices of a rectified
5-simplex; i.e., each vertex being at the center of an edge of a
regular 5-simplex.
A MATLAB program to generate these rules ({\tt c5\_22\_4.m}) is included in the supplemental materials.

\begin{table}
\caption{22 Point Rules of Degree 4 in 5 Dimensions.}\label{tab:5_22_4_p}
\begin{center}
\begin{tabular}{r@{}cr@{}cr@{}cr@{}cr@{}ccccc}
 $ x_{1}$& &$ x_{2}$& &$ x_{3}$& &$ x_{4}$& &$ x_{5}$& &Weight& \# Points\\
\hline
    0  & &   0  & &   0  & &   0  & &   0  &     &$ W_0$& 1 \\
 $  c $& &$  c $& &$  c $& &$  c $& &$  c $&     &$ W_1$& 1 \\
 $(-h $& &$  a $& &$  a $& &$  a $& &$  a $&$)_S$&$ W_1$& 5 \\
 $(-b $& &$ -b $& &$ -b $& &$  g $& &$  g $&$)_S$&$ W_2$& 10 \\
 $( e $& &$ -f $& &$ -f $& &$ -f $& &$ -f $&$)_S$&$ W_2$& 5 \\
\end{tabular}

\begin{tabular}{clll}\\
        &\hspace{3.5em}$E_5^{r^2}\T $    &  \hspace{3.5em}$E_5^r    $      &  \hspace{3.5em}$S_5$         \\
\hline
$  a \T$&$ (2\sqrt{3}-\phantom{2}\sqrt{2})/10 $&$ (\phantom{1} 4\sqrt{3} - 2\sqrt{2})/5 $&$ (\phantom{2} \sqrt{6} - 1)/15 $\\  
$  b   $&$ (2\sqrt{3}-\phantom{2}\sqrt{2})/5  $&$ (\phantom{1} 8\sqrt{3} - 4\sqrt{2})/5 $&$ (2\sqrt{6} - 2)/15  $\\  
$  c   $&$ \sqrt{1/2}                         $&$ \sqrt{8}                              $&$ 1/3                 $\\  
$  e   $&$ (4\sqrt{3}-2\sqrt{2})/5            $&$ (16\sqrt{3} - 8\sqrt{2})/5            $&$ (4\sqrt{6} - 4)/15  $\\  
$  f   $&$ (\phantom{2}\sqrt{3}+2\sqrt{2})/5  $&$ (\phantom{1} 4\sqrt{3} + 8\sqrt{2})/5 $&$ (\phantom{2} \sqrt{6} + 4)/15 $\\  
$  g   $&$ (3\sqrt{3}+\phantom{2}\sqrt{2})/5  $&$ (12\sqrt{3} + 4\sqrt{2})/5            $&$ (3\sqrt{6} + 2)/15  $\\
$  h   $&$ (8\sqrt{3}+\phantom{2}\sqrt{2})/10 $&$ (16\sqrt{3} + 2\sqrt{2})/5            $&$ (4\sqrt{6} + 1)/15  $\\
$  W_0 $&$ \pi^{5/2}/4                         $&$  28\pi^2                              $&$ 2\pi^2/105          $\\
$  W_1 $&$ \pi^{5/2}/18                        $&$  8\pi^2/3                             $&$ 4\pi^2/105          $\\
$  W_2 $&$ \pi^{5/2}/36                        $&$  4\pi^2/3                             $&$ 2\pi^2/105          $\\
$  r_1 $&$ \sqrt{5/2}                         $&$  \sqrt{40}                            $&$ \sqrt{5}/3          $\\
$  r_2 $&$ 2                                  $&$    8                                  $&$ \sqrt{8}/3          $\\
\end{tabular}
\end{center}
\end{table}

\subsection{Degree 4, dimension 6, 28 point rules}

\begin{figure}
  \centering
  \includegraphics[width=3in]{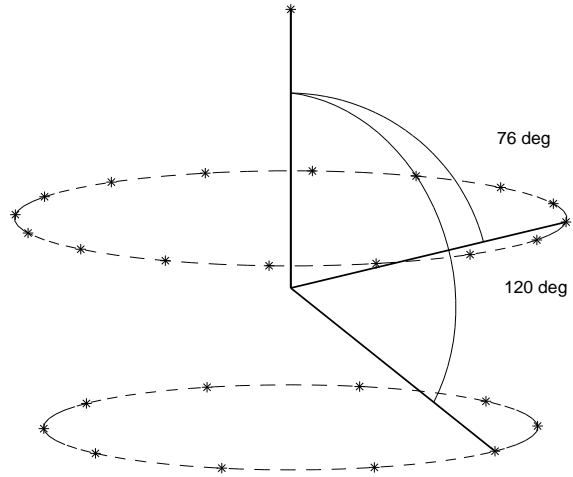}
  \caption{Configuration of the 27 non-central points for the rules of degree 4 in 6 dimensions.  
}
  \label{fig:6_28_4}
\end{figure}

Each of these new formulas has a central point and 27 points
all at the same radius.
The rule is shown in Table \ref{tab:6_28_4}.
Each point on the shell has 16 near neighbors
(~76 degrees away) and 10 more distant neighbors (120 degrees).
This configuration is suggested in Figure \ref{fig:6_28_4}, which
shows the points in terms of their angular distance from a
chosen point, though of course their distances from each other cannot
be shown realistically.

\begin{table}
\caption{28 Point Rules of Degree 4 in 6 Dimensions.}\label{tab:6_28_4}
\begin{center}
\begin{tabular}{r@{}cr@{}cr@{}cr@{}cr@{}cr@{}cccc}
$ x_{1}$& &$ x_{2}$& &$ x_{3}$& &$ x_{4}$& &$ x_{5}$& &$ x_{6}$& &Weight   & Radius & \# Points \\
\hline

    0  & &   0  & &   0  & &   0  & &   0  & &   0  &     &$ W_0$ &0& 1 \\
 $ -c $& &$\pm e $& &   0  & &   0  & &   0  & &   0  &     &$ W_1$ &$r  $& 2 \\
 $ -c $& &   0  & &$\pm(b $& &$  b $& &$  b $& &$ -b $&$)_S$&$ W_1$ &$r  $& 8 \\
 $  a $& &$ -b $& &$\pm(b $& &$  b $& &$  b $& &$  b $&)    &$ W_1$ &$r  $& 2 \\
 $  a $& &$ -b $& &$ (b $& &$  b $& &$ -b $& &$ -b $&$)_S$&$ W_1$ &$r  $& 6 \\
 $  f $& &   0  & &   0  & &   0  & &   0  & &   0  &     &$ W_1$ &$r  $& 1 \\
 $  a $& &$  b $& &$\pm(e $& &   0  & &   0  & &   0  &$)_S$&$ W_1$ &$r  $& 8 \\
\end{tabular}

\begin{tabular}{lcccc}\\
          &  $ E_6^{r^2}\T $  &  $   E_6^r  $  &  $S_6       $  \\
\hline
$  a     $& $         1/2$ & $  \sqrt{9/2}$\T & $ \sqrt{1/20}$\\ 
$  b     $& $  \sqrt{3/4}$ & $ \sqrt{27/2}$ & $ \sqrt{3/20}$\\ 
$  c     $& $           1$ & $   \sqrt{18}$ & $  \sqrt{1/5}$\\ 
$  e     $& $    \sqrt{3}$ & $   \sqrt{54}$ & $  \sqrt{3/5}$\\ 
$ f=r    $& $           2$ & $   \sqrt{72}$ & $  \sqrt{4/5}$\\ 
$  W_0   $& $     \pi^3/4$ & $     50\pi^3$ & $    \pi^3/96$\\ 
$  W_1   $& $    \pi^3/36$ & $  70\pi^3/27$ & $  5\pi^3/864$\\ 
\end{tabular}
\end{center}
\end{table}

\subsection{Degree 4, dimension 7, 38 point rules}

\begin{figure}
  \centering
  \includegraphics[width=3in]{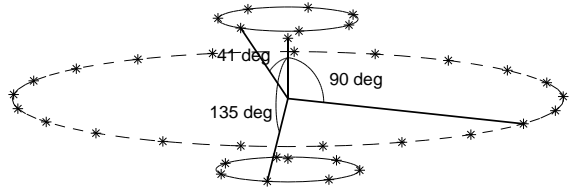}
  \caption{Configuration of the 37 non-central points for the rules of degree 4 in 7 dimensions.}
  \label{fig:hoops-1-1-7-7-21}
\end{figure}

We found 38 point rules for $E_7^{r^2}$ and $S_7$, each with
two negative weights. A standard measure of the stability of an
integration rule is the sum of the absolute value of the weights,
divided by the sum of the weights, which is a worst-case round-off error
magnification factor \cite{genz-keister96interpolatory}.  These rules
have stability factors of 7.18 for $E_7^{r^2}$ and 8.55 for $S_7$.
They have some symmetry, with the center point, a centered
shell of 21 points, and two offset irregular 7-simplices.
The configuration of the points with respect to one of the negative
weight points is suggested in Figure \ref{fig:hoops-1-1-7-7-21}.
The rules are
shown in Table \ref{tab:7_38_4}.

\begin{table}
\caption{38 Point Rules of Degree 4 in 7 Dimensions.}\label{tab:7_38_4}
\begin{center}\begin{tabular}{r@{}cr@{}cr@{}cr@{}cr@{}cr@{}cr@{}cccc}
$x_{1}$& &$x_{2}$& &$x_{3}$& &$x_{4}$& &$x_{5}$& &$x_{6}$& &$x_{7}$& &Weight&\# Points \\
\hline
   0  & &   0  & &   0  & &   0  & &   0  & &   0  & &   0  & &$ \phantom{-}W_0  $& 1  \\
$  c $& &$  c $& &$  c $& &$  c $& &$  c $& &$  c $& &$  c $& &$ -W_2  $& 1  \\
$ -b $& &$ -b $& &$ -b $& &$ -b $& &$ -b $& &$ -b $& &$ -b $& &$ -W_1  $& 1  \\
$(\hspace{1ex}f $& &$ -e $& &$ -e $& &$ -e $& &$ -e $& &$ -e $& &$ -e $&$)_S$&$ \phantom{-}W_3  $& 7  \\
$ (\hspace{1ex}h $& &$  a $& &$  a $& &$  a $& &$  a $& &$  a $& &$  a $&$)_S$&$ \phantom{-}W_4  $& 7  \\
$(-i $& &$ -i $& &$  g $& &$  g $& &$  g $& &$  g $& &$  g $&$)_S$ &$ \phantom{-}W_5  $& 21\\
\end{tabular}

\begin{tabular}{crc}\\   
        &  $   E_7^{r^2}\T\hspace{2.5em} $ &  $   S_7       $   \\ 
\hline
$  a   $&  0.2286166663871 &  0.0974824740891  \\ 
$  b   $&  0.2590817563916 &  0.1104728321147  \\ 
$  c   $&  0.3117777721419 &  0.1329424887288  \\ 
$  e   $&  0.4422503418055 &  0.1885761793629  \\ 
$  f   $&  0.4505846393780 &  0.1921299357884  \\ 
$  g   $&  0.7531484451994 &  0.3211435760773  \\ 
$  h   $&  1.0981884332902 &  0.4682691213418  \\ 
$  i   $&  1.8927504201541 &  0.8070714909185  \\ 
$  W_0 $& 59.8014451908073 &  5.2337832579847  \\ 
$  W_1 $& 89.9014937680773 &  9.4465413692728  \\ 
$  W_2 $& 79.9432767398149 &  8.4001659957515  \\ 
$  W_3 $& 11.6616239025637 &  1.2253635397056  \\ 
$  W_4 $& 11.0688850060780 &  1.1630805645052  \\ 
$  W_5 $&  0.2803313076587 &  0.0294562546617  \\   
\end{tabular}
\end{center}
\end{table}

\section{Degree 5 rules}

\subsection{Degree 5, dimension 4, 23 point rules}

This new family of degree 5 rules 
is shown in Table \ref{tab:4_23_5}.
Other than the central point, all points are at the same radius.
However, two of those points have lower weight than the others.
The rule for $E_5^r$ has fewer points than for the known rules.
A 22 point rule of degree 5 was known for $S_n$
\cite[$S_n$:5-1]{stroud71approximate}.
\begin{table}
\caption{23 Point Rules of Degree 5 in 4 Dimensions.}\label{tab:4_23_5}
\begin{center}\begin{tabular}{r@{}cr@{}cr@{}cr@{}ccccc}
$x_{1}$& &$x_{2}$& &$x_{3}$& &$x_{4}  $& &  Weight &Radius&\# Points \\
\hline
    0  & &   0   & &   0  & &   0     & &$W_{0}  $&0& 1  \\
$\pm h$& &   0   & &   0  & &   0     & &$W_{2}  $&$r$& 2  \\
    0  & &$\pm h$& &   0  & &   0     & &$W_{1}  $&$r$& 2  \\
                                                  
$\pm c$&&$\pm( b$& &$ -a $&) &$\pm c $& &$W_{1}  $&$r$& 8  \\
$\pm c$& &$ \pm( b $& &$ e$&) &  0   & &$  W_{1}  $&$r$& 4  \\
   0   & &$\pm(a$& &$ -g $&) &   0   & &$  W_{1}  $&$r$& 2  \\
   0   & &$\pm(a$& &$  b $&) &$\pm f$& &$  W_{1}  $&$r$& 4  \\
\end{tabular}
\end{center}
\begin{center}
\begin{tabular}{lcccc}\\
         &  $   E_4^{r^2}\T $ &  $ E_4^r        $ &  $   S_4         $  \\ 
\hline
 $  a   $&$ \sqrt{1/3}\T $   &$ \sqrt{14/3}    $&$  \sqrt{1/12}    $\\          
 $  b   $&$ \sqrt{2/3} $     &$ \sqrt{28/3}    $&$  \sqrt{1/6}     $\\          
 $  c   $&$ 1          $     &$ \sqrt{14}      $&$  \sqrt{1/4}     $\\          
 $  e   $&$ \sqrt{4/3} $     &$ \sqrt{56/3}    $&$  \sqrt{1/3}     $\\          
 $  f   $&$ \sqrt{2  } $     &$ \sqrt{28}      $&$  \sqrt{1/2}     $\\          
 $  g   $&$ \sqrt{8/3} $     &$ \sqrt{112/3}   $&$  \sqrt{2/3}     $\\          
 $  h=r   $&$ \sqrt{3  } $     &$ \sqrt{42}      $&$  \sqrt{3/4}     $\\          
 $  W_0 $&$ \pi^2/3    $     &$ 44\pi^2/7      $&$  \pi^2/18       $\\          
 $  W_1 $&$ \pi^2/32   $     &$ 15\pi^2/56     $&$  \pi^2/48       $\\          
 $  W_2 $&$ \pi^2/48   $     &$ 5\pi^2/28      $&$  \pi^2/72       $\\            
\end{tabular}
\end{center}
\end{table}

\subsection{Degree 5, dimension 6, 44 point rule}
This new rule has points supported by two spheres,  
as shown in Table \ref{tab:6_44_5}.
\begin{table}
\caption{44 Point Rule of Degree 5 for  $E_6^r$.}\label{tab:6_44_5}
\begin{center}\begin{tabular}{r@{}cr@{}cr@{}cr@{}cr@{}cr@{}cccc}
$x_{1}$& &$x_{2}$& &$x_{3}$& &$x_{4}$& &$x_{5}$& &$x_{6}$& &Weight&Radius&\# Points \\
\hline
  (0  & &   0  & &   0  & &   0  & &   0  & &$\pm b$&$)_S$&$ \phantom{-}W_{1}  $&$r_1$& 12  \\
$  (a $& &$  a $& &$  a $& &$  a $& &$  a $& &$ -a $&$)_S$&$ \phantom{-}W_{2}  $&$r_2$& 6  \\
$ (a $& &$  a $& &$  a $& &$ -a $& &$ -a $& &$ -a $&$)_S$&$ \phantom{-}W_{2}  $&$r_2$& 20  \\
$  (a $& &$ -a $& &$ -a $& &$ -a $& &$ -a $& &$ -a $&$)_S$&$ \phantom{-}W_{2}  $&$r_2$& 6  \\
\end{tabular}
\end{center}
\begin{center}
\begin{tabular}{cc}\\
         &  $   E_6^r     $   \\
\hline
$  a    $&  4.84099298434420  \\
$  b=r_1$&  5.40578920173885  \\
$   W_1 $&  274.495347525855  \\
$   W_2 $&  13.3377822289287  \\
$    r_2$&  11.8579626600364\\
\end{tabular}
\end{center}
\end{table}

\section{Degree 6 rules}
\subsection{Degree 6, dimension 2, 10 point rule}

\begin{figure}
  \centering
 
  \includegraphics[width=5in,angle=0]{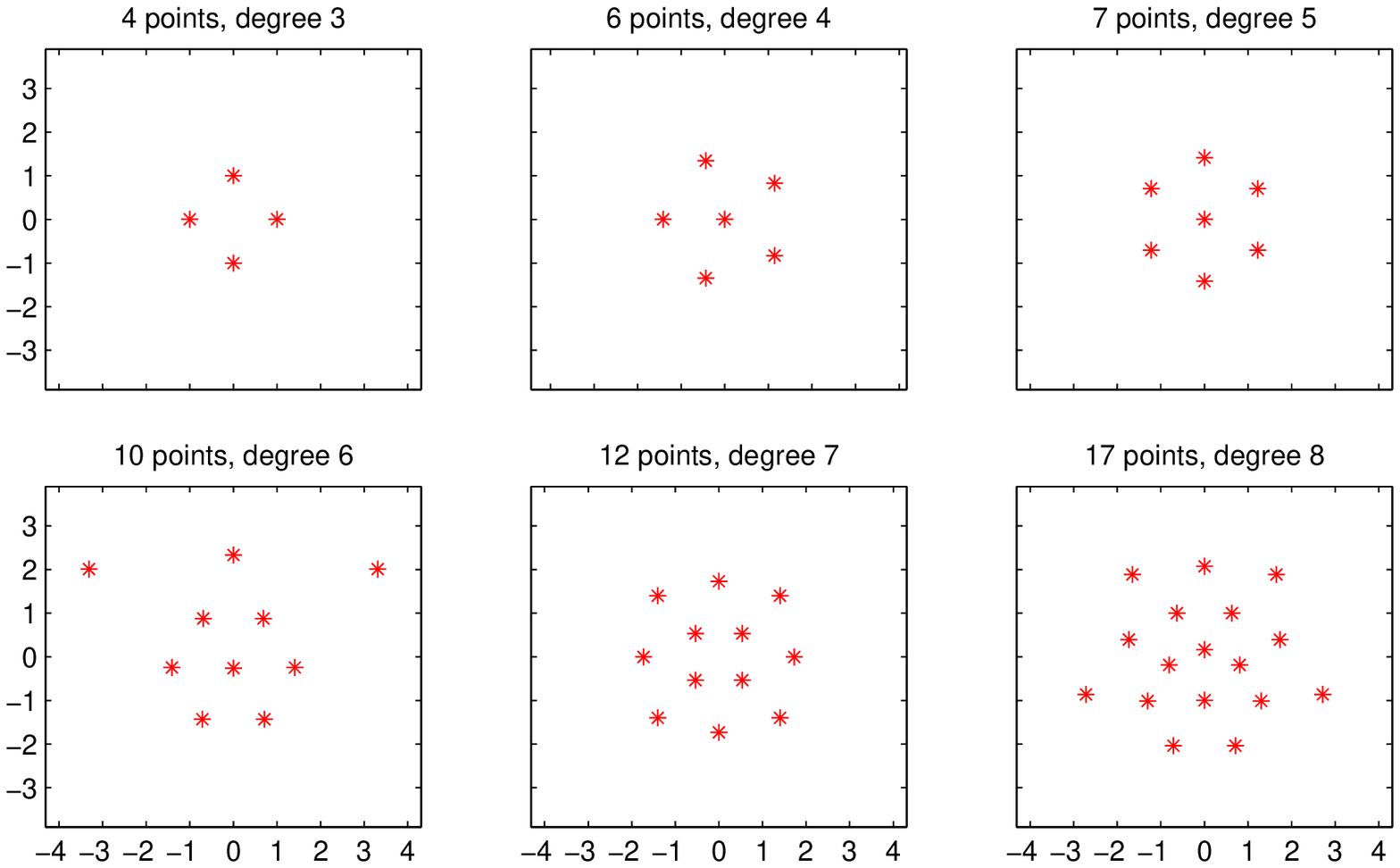}
  \caption{Points for $E_2^{r^2}$ rules.  The rule of degree 8 is new.}
  \label{fig:plot3e2}
\end{figure}

A rule was found for $E_2^{r^2}$ with 10 points, achieving Stroud's lower bound \eqref{eq:stroud}.
This rule was known, but unpublished\cite{cools03cubature,jankewitz98}.
The points and weights are shown in Table \ref{tab:2_10_6}.
The point layout is similar to that in
the 10 point rule for $S_2$ by Wissmann and Becker
\cite[$S_2$:6-1]{wissmann-becker86cubature}.  
The points are shown in Figure \ref{fig:plot3e2}, along with those for
known formulas of degree 3, 4, 5, and 7, and the new formula of degree
8 discussed below.
Note in the figure that the rules of odd degree have
central symmetry (for every point $x$, there is also
a point $-x$ with the same weight), while those of even degree are
only bilaterally symmetric.

\begin{table}
\begin{center}
\caption{10 Point Rule of Degree 6 for $E_2^{r^2}$.}\label{tab:2_10_6}
\begin{tabular}{rrcc}
$ x_{1}$\hspace{3em} &$ x_{2}$\hspace{3em} &Weight&Radius   \\
\hline
$\pm$3.314013565941806 &  2.014171295633760 & 0.000757833922865 &   3.87809 \\
$\pm$1.411670545911536 & -0.242569904073576 & 0.236161927729435 &   1.43236 \\
$\pm$0.713033732783175 & -1.432390280414699 & 0.146082553662775 &   1.60005 \\
$\pm$0.691608815107559 &  0.877693534044218 & 0.485399260031153 &   1.11744 \\
     0.000000000000000 & -0.261367769356158 & 1.387418367858287 &   0.26137 \\
     0.000000000000000 &  2.335832264987514 & 0.017371135039050 &   2.33583\\

\end{tabular}
\end{center}
\end{table}

\subsection{Degree 6, dimension 2, 11 point rule}

This new rule for $E_2^r$ has 11 points, with bilateral
symmetry.  The points and weights are shown in Table \ref{tab:2_11_6E1}.
This rule come close to pentagonal symmetry, but we were unable to
adjust it for pentagonal symmetry.

\begin{table}
\begin{center}
\caption{11 Point Rule of Degree 6 for $E_2^r$.}\label{tab:2_11_6E1}
\begin{tabular}{rrcr}
$ x_{1}$\hspace{3em} &$ x_{2}$\hspace{3em} &Weight&Radius   \\
\hline
      0.000000000000000 &  0.000000000000000 & 3.927702275194840 &  0.00000 \\
      0.000000000000000 & 10.299713185154499 & 0.003846684331349 & 10.29971 \\
      0.000000000000000 & -3.895765525253948 & 0.474246212300936 &  3.89577 \\
$\pm$10.311630315898372 &  3.397224688449697 & 0.002841012046587 & 10.85683 \\
$\pm$ 6.251012172182811 & -8.794364006109971 & 0.002944454683352 & 10.78962 \\
$\pm$ 3.752487980256190 & -1.228482827331175 & 0.460111970539923 &  3.94846 \\
$\pm$ 2.312667676618243 &  3.141828043257887 & 0.472797630406369 &  3.90122\\
\end{tabular}
\end{center}
\end{table}

\section{Degree 7 rules}

We found four rules of degree 7 with fewer points than previously reported.

\subsection{Degree 7, dimension 6, 127 point rules}

Each of these new rules has a central point and two shells.  The inner shell
has 54 points.  Each of those has 10 nearest neighbors in that shell
(60 degrees away), and 16, 16, 10, and 1, successively further away.
The outer shell has 72 points.  Each of those has 20 nearest
neighbors, and 30, 20, and 1, successively further away.  
The configuration of points is shown in Table \ref{tab:6_127_7}.
\begin{table}\caption{127 Point Rules of Degree 7 in 6 Dimensions.}\label{tab:6_127_7}
\begin{center}\begin{tabular}{r@{}cr@{}cr@{}cr@{}cr@{}cr@{}ccccc}
$ x_{1}$& &$ x_{2}$& &$ x_{3}$& &$ x_{4}$& &$ x_{5}$& &$ x_{6}$& &Weight&Radius&   \# Points \\
\hline
  0     & &  0     & &  0     & &  0     & &  0     & &  0     & &$ W_0$&0& 1  \\
$\pm g    $& &  0     & &  0     & &  0     & &  0     & &  0     &&$ W_1$&$r_1$&  2\\
$\pm c    $&&($\pm f    $& &  0     & &  0     & &  0     & &  0     &$)_S$&$ W_1$&$r_1$&  20 \\
$\pm(a    $& &$ b    $& &$ b    $& &$ b    $& &$ b    $& &$ b    $&)&$ W_1$&$r_1$&  2 \\
$\pm(a    $&&$( b    $& &$ b    $& &$ b    $& &$-b    $& &$-b    $&$)_S)$&$ W_1$&$r_1$&  20 \\
$\pm(a    $&&$( b    $& &$-b    $& &$-b    $& &$-b    $& &$-b    $&$)_S)$&$ W_1$&$r_1$&  10 \\
$\pm(h    $&&$( e    $& &$ e    $& &$ e    $& &$ e    $& &$-e    $&$)_S)$&$ W_2$&$r_2$& 10 \\
$\pm(h    $&&$( e    $& &$ e    $& &$-e    $& &$-e    $& &$-e    $&$)_S)$&$ W_2$&$r_2$& 20 \\
$\pm(h    $& &$-e    $& &$-e    $& &$-e    $& &$-e    $& &$-e    $&)&$ W_2$&$r_2$& 2 \\
 0   & &$( \pm i    $& &$\pm i    $& &  0     & &  0     & &  0     &$)_S$&$ W_2$&$r_2$& 40 \\
\end{tabular}

\begin{tabular}{llc}\\
        &\hspace{3em} $   E_6^{r^2}\T $    &  $   S_6   $  \\
\hline
$g=r_1 $& $\sqrt{(4-\sqrt{6})\times2}\TT$ & $\sqrt{2/3}$  \\ 
$  c   $& $\sqrt{(4-\sqrt{6})/2}       $ & $\sqrt{1/6}$  \\ 
$  f   $& $\sqrt{(4-\sqrt{6})\times3/2}$ & $\sqrt{1/2}$  \\ 
$  a   $& $\sqrt{(4-\sqrt{6})/8}       $ & $\sqrt{1/24}$ \\ 
$  b   $& $\sqrt{(4-\sqrt{6})\times3/8}$ & $\sqrt{1/8}$  \\ 
$  e   $& $\sqrt{(6+\sqrt{6})/8}       $ & $\sqrt{1/8}$  \\ 
$  h   $& $\sqrt{(6+\sqrt{6})\times3/8}$ & $\sqrt{3/8}$  \\ 
$  i   $& $\sqrt{(6+\sqrt{6})/2}       $ & $\sqrt{1/2}$  \\ 
$r_2   $& $\sqrt{6+\sqrt{6}}           $ &        1      \\ 
                                                            
$  W_0 $& $(16-\sqrt{6})\pi^3/100      $ & $\pi^3/240$  \\  
$  W_1 $& $(68+27\,\sqrt{6})\pi^3/9000 $ & $\pi^3/480$  \\  
$  W_2 $& $(54-19\,\sqrt{6})\pi^3/9000 $ & $\pi^3/1440$ \\    
\end{tabular}
\end{center}
\end{table}

\subsection{Degree 7, dimension 7, 183 point rules}\label{r7_183_7}
In this case, we initialized a search with 226 points, and this new rule
was found---the weights on the remaining 43 points having been driven
to zero.  It does not quite attain M{\"o}ller's lower bound of
$N=n/3(n^2+3n+8)=182$ for a degree seven formula
\cite{moller79nodes,lu-darmofal04gaussian}.

The rule has a central point, one shell of 56 points, and a second
shell of 126 points.  The inner shell is laid out the same as for the
57 point formula of degree 5 by Stroud \cite[$E_n^{r^2}$:5-1]{stroud71approximate}.
The points on the outer shell have vertex symmetry, but we have been
unable to relate them to a known polytope.

The points are shown in Table \ref{tab:7_183_7}.
We found closed form expressions for the points on the outer
shell and for the radius $r_1$ of the inner shell directly from their
simple continued fractions.
We were then able to find expressions for the ratios of the
remaining coordinates to $r_1$.
Maxima was then able to solve for the coordinates using the expressions
for the points and three of the moment constraint equations.
\begin{table}
\caption{183 Point Rules of Degree 7 in 7 Dimensions.}\label{tab:7_183_7}
\ifels
\vspace*{-15pt}
\fi
\begin{center}\begin{tabular}{r@{}cr@{}cr@{}cr@{}cr@{}cr@{}cr@{}lccc}
$ x_{1}$& &$ x_{2}$& &$ x_{3}$& &$ x_{4}$& &$ x_{5}$& &$ x_{6}$& &$ x_{7}$& &Weight&Radius&\# Points\\
\hline
   0  & &   0  & &   0  & &   0  & &   0  & &   0  & &   0  & &$ W_0 $&0& 1 \\
$\pm(-m $& &   0  & &   0  & &   0  & &   0  & &   0  & &   0  &)&$ W_1 $&$r_1$&    2 \\
$\pm(-c $& &$  k $& &   0  & &   0  & &   0  & &   0  & &   0  &)&$ W_1 $&$r_1$&    2 \\
$\pm(-c $& &$ -f $&&$(\pm i $& &   0  & &   0  & &   0  & &   0  &$)_S)$&$ W_1 $&$r_1$&    20 \\
$\pm(-c $& &$  a $& &$  e $& &$  e $& &$  e $& &$  e $& &$  e $&)&$ W_1 $&$r_1$&    2 \\
$\pm(-c $& &$  a $&&$(  e $& &$  e $& &$  e $& &$ -e $& &$ -e $&$)_S)$&$ W_1 $&$r_1$&    20 \\

$\pm(-c $& &$  a $&&$(  e $& &$ -e $& &$ -e $& &$ -e $& &$ -e $&$)_S)$&$ W_1 $&$r_1$&    10 \\
$\pm(j $& &$  p $& &   0  & &   0  & &   0  & &   0  & &   0  &)&$ W_2 $&$r_2$&   2 \\
$\pm(j $& &$  b $& &$  g $& &$  g $& &$  g $& &$  g $& &$  g $&)&$ W_2 $&$r_2$&   2 \\
$\pm(j $& &$  b $&&$ ( g $& &$  g $& &$  g $& &$ -g $& &$ -g $&$)_S)$&$ W_2 $&$r_2$&   20 \\
$\pm(j $& &$  b $&&$(  g $& &$ -g $& &$ -g $& &$ -g $& &$ -g $&$)_S)$&$ W_2 $&$r_2$&   10 \\
$\pm(j $& &$ -h $&&$(\pm o $& &   0  & &   0  & &   0  & &   0  &$)_S)$&$ W_2 $&$r_2$&   20 \\
$\pm(0$  & &$  l $&&$ ( g $& &$  g $& &$  g $& &$  g $& &$ -g $&$)_S)$&$ W_2 $&$r_2$&   10 \\
$\pm$(0  & &$  l $&&$ ( g $& &$  g $& &$ -g $& &$ -g $& &$ -g $&$)_S)$&$ W_2 $&$r_2$&   20 \\
$\pm($0  & &$  l $& &$ -g $& &$ -g $& &$ -g $& &$ -g $& &$ -g $&)&$ W_2 $&$r_2$&   2 \\
      0  & &   0  &&$(\pm o $& &$\pm o $& &   0  & &   0  & &   0  &$)_S$&$ W_2 $&$r_2$&   40 \\
\end{tabular}
  \begin{tabular}{rlll}\\
        & \hspace{3em} $   E_7^{r^2}  \T  $    & \hspace{3em} $   S_7       $                        \\ 
\hline
n$ m=r_1$& $\sqrt{(9-4\,\sqrt{3})\times3/2}$\TT&$\sqrt{(117-4\,\sqrt{78})\times3/377}$               \\  
$  c   $& $\sqrt{(9-4\,\sqrt{3})/6}$          &$\sqrt{(117-4\,\sqrt{78})/1131}$                     \\  
$  k   $& $\sqrt{(9-4\,\sqrt{3})\times4/3}$   &$\sqrt{(117-4\,\sqrt{78})\times8/1131}$              \\  
$  f   $& $\sqrt{(9-4\,\sqrt{3})/3}$          &$\sqrt{(117-4\,\sqrt{78})\times2/1131}$              \\  
$  i   $& $\sqrt{9-4\,\sqrt{3}}$              &$\sqrt{(117-4\,\sqrt{78})\times2/377}$               \\  
$  a   $& $\sqrt{(9-4\,\sqrt{3})/12}$         &$\sqrt{(117-4\,\sqrt{78})/2262}$                     \\  
$  e   $& $\sqrt{(9-4\,\sqrt{3})/4}$          &$\sqrt{(117-4\,\sqrt{78})/754}$                      \\  
$  j   $& $\sqrt{(\sqrt{3}+6)/3}$             &$\sqrt{(\sqrt{78}+78)/273}$                          \\  
$  p   $& $\sqrt{(\sqrt{3}+6)\times2/3}$      &$\sqrt{(\sqrt{78}+78)\times2/273}$                   \\  
$  b   $& $\sqrt{(\sqrt{3}+6)/24}$            &$\sqrt{(\sqrt{78}+78)/2184}$                         \\  
$  g   $& $\sqrt{(\sqrt{3}+6)/8}$             &$\sqrt{(\sqrt{78}+78)/728}$                          \\  
$  h   $& $\sqrt{(\sqrt{3}+6)/6}$             &$\sqrt{(\sqrt{78}+78)/546}$                          \\  
$  o   $& $\sqrt{(\sqrt{3}+6)/2}$             &$\sqrt{(\sqrt{78}+78)/182}$                          \\  
$  l   $& $\sqrt{(\sqrt{3}+6)\times3/8}$      &$\sqrt{(\sqrt{78}+78)\times3/728}$                   \\  
$  r_2 $& $\sqrt{\sqrt{3}+6}$                 &$\sqrt{(\sqrt{78}+78)/91}$                           \\  
$  W_0 $&$(144-35\,\sqrt{3}) \,\pi^{7/2}/1089 $&$(6912-7\times2^{11/2}\,\sqrt{39})\,\pi^3/2264031       $\\
$  W_1 $&$(675+388\,\sqrt{3})\,\pi^{7/2}/95832$&$(104598+1085\times2^{7/2}\,\sqrt{39})\,\pi^3/124521705 $\\
$  W_2 $&$(90-37\,\sqrt{3})  \,\pi^{7/2}/23958$&$(101088-235\times2^{9/2}\,\sqrt{39})\,\pi^3/124521705  $ \\

\end{tabular}
\end{center}
\end{table}

\section{Degree 8 rules}

\subsection{Degree 8, dimension 2, 17 point rule}

We found 17 point rules of degree 8 for all three integrals with all positive weights and
bilateral symmetry.
For details, see the supplemental material.
A 16 point rule of degree 8 for $S_3$ was found by
Wissmann and Becker
\cite{wissmann-becker86cubature}.  We were unable, even using variations of
that rule as starting guesses, to find a similar rule for $E_2^r$ or $E_2^{r^2}$.
\section{Degree 9 rules}\label{deg9}

\subsection{Degree 9, dimension 4, 124 point rule}

We found a 124 point rule for $E_4^{r^2}$ with negative weights (stability
factor 15.4) and central symmetry, but no central point.  
We also found a 125 point rule for the same integrals with central
symmetry and a central
point.  It also has negative weights, but a somewhat better stability
factor of 8.1.
For details, see the supplemental material.

\section{Summary}\label{summary}
\subsection{Listings}
The new rules are listed in Tables \ref{tab:summaryE2},
\ref{tab:summaryE1}, and \ref{tab:summaryS}.
In addition to those described above, we found many rules with only
bilateral symmetry or no apparent symmetry, the details for which
appear only in the supplemental material.
Symmetry of ``$x_2,x_3$" indicates a rule closed under sign changes in
both of the indicated coordinates. 
Rules with the symmetry of a known polytope are indicated by that
polytope.
``Vertex'' indicates symmetry with respect to the exchange of any two
noncentral points, but that the polytope has not been identified.

The ``Quality'' of a rule is given using the notation introduced in
\cite{lyness-jespersen75triangle}.  The first letter is 
P if all weights are positive,
or N if some weights are negative.
For integral $S_n$, there is a second letter,
which is I if all points are inside the region, or B if some are on
the boundary, or O if some points are outside the region.

Also shown is the M{\"o}ller Lower Bound (MLB) for the number of
points in a rule of the given degree from \eqref{eq:mlb}, and the
smallest known rule of the given degree or the next
higher degree.  
The new rules with points supported by one or two spherical shells are
very efficient---within three points of the M{\"o}ller lower bound.
Those with little or no symmetry are much less efficient, with over 40
percent more points than the M{\"o}ller lower bound in the median; though still better than the previously known rules, with the exceptions noted in the tables.

In most odd degree formulas, points are supported by a few spherical shells,
with all weights positive.  
Most even degree formulas lack symmetry, and they had more negative weights.
We were unable to find rules for $E_n^r$ and sometimes even $S_n$ corresponding to some of the
rules for $E_n^{r^2}$.
\def\clap#1{\hbox to 0pt{\hss\;\;#1\hss}}
\begin{table}
  \caption{ 25 New Cubature Rules for $E_n^{r^2}$.}\label{tab:summaryE2}
  \begin{center}
    \begin{tabular}{ccc @{\hspace{1ex}}c@{\hspace{-4pt}} c@{\hspace{2ex}}ccccc}
     
      \multicolumn{6}{c}{New Rule} &
      \multirow{2}{*}{$N_{MLB}$}&
      \multicolumn{3}{c} {Smallest Previous Rule} \\
      \cmidrule(lr){1-6}  \cmidrule(lr){8-10}
      $n$ & $N$ & $d$ & Shells & Quality & Symmetry & & $N$ & $d$ & Source \\
      \cmidrule(lr){1-6} \cmidrule(lr){7-7}  \cmidrule(lr){8-10}
        2 &     17 &  8 &           &      P &   bilateral & 15 &   18 & 9 & \cite{haegemans-piessens77seven}  \\
  2 &     24 & 10 &           &      N &   bilateral & 21 &   25 &11 & \cite{haegemans-piessens76eleven}  \\
  3 &     10 &  4 &           &      P &   $x_2,x_3$ & 10 &   13 & 5 & \cite[VII]{stroud-secrest63spherically}  \\
  3 &     22 &  6 &           &      P &   bilateral & 20 &   27 & 7 & \cite[$E_n^{r^2}$:7-1]{stroud71approximate}  \\
  3 &    220 & 14 &           &      N &        none &120 &  288 &14 & \cite[$E_3^{r^2}$:14-1]{stroud71approximate}  \\
  3 &    234 & 15 &           &      N &        none &140 &\multicolumn{3}{c}{none} \\
  4 &     16 &  4 &    1+10+5 &      P &   4-simplex & 15 &   22 & 5 & \cite[$E_n^{r^2}$:5-1]{stroud71approximate}  \\
  4 &     16 &  4 &     1+6+9 &      P &   $x_1,x_3$ & 15 &   22 & 5 & \cite[$E_n^{r^2}$:5-1]{stroud71approximate}  \\
  4 &    23\clap{*} &  5 &      1+22 &      P &      vertex & 21 &   22 & 5 & \cite[$E_n^{r^2}$:5-1]{stroud71approximate}  \\
  4 &     43 &  6 &           &      P &   bilateral & 35 &   49 & 7 & \cite[$E_n^{r^2}$:7-1]{stroud71approximate}  \\
  4 &    105 &  8 &           &      N &        none & 70 &  193 & 9 & \cite[$E_n^{r^2}$:9-1]{stroud71approximate}  \\
  4 &    124 &  9 &           &      N &     central & 91 &  193 & 9 & \cite[$E_n^{r^2}$:9-1]{stroud71approximate}  \\
  4 &    125 &  9 &           &      N &     central & 91 &  193 & 9 & \cite[$E_n^{r^2}$:9-1]{stroud71approximate}  \\
  4 &    213 & 10 &           &      N &        none &126 &  417 &11 & \cite[$E_n^{r^2}$:11-1]{stroud71approximate}  \\
  5 &     22 &  4 &    1+6+15 &      P &   5-simplex & 21 &   32 & 5 & \cite[$E_n^{r^2}$:5-1]{stroud71approximate}  \\
  5 &     80 &  6 &           &      P &        none & 56 &   83 & 7 & \cite[$E_n^{r^2}$:7-1]{stroud71approximate}  \\
  5 &    224 &  8 &           &      N &        none &126 &  395 & 9 & \cite[CUT8]{adurthi-singla-singh12cut}  \\
  6 &     28 &  4 &      1+27 &      P &      vertex & 28 &   44 & 5 & \cite[$E_n^{r^2}$:5-1]{stroud71approximate}  \\
  6 &    127 &  7 &   1+54+72 &      P &     central &124 &  137 & 7 & \cite[$E_n^{r^2}$:7-1]{stroud71approximate}  \\
  7 &     38 &  4 &  1+8+8+21 &      N &    see text & 36 &   57 & 5 & \cite[$E_n^{r^2}$:5-1]{stroud71approximate}  \\
  7 &    183 &  7 &  1+56+126 &      P &     central &182 &  227 & 7 & \cite[$E_n^{r^2}$:7-1]{stroud71approximate}  \\
  8 &    339 &  6 &           &      N &        none &165 &  705 & 7 & \cite[$E_n^{r^2}$:7-3]{stroud71approximate}  \\
  9 &     76 &  4 &           &      P &        none & 55 &  111 & 5 & \cite[I]{lu-darmofal04gaussian}  \\
 10 &     96 &  4 &           &      P &        none & 66 &  133 & 5 & \cite[I]{lu-darmofal04gaussian}  \\
 11 &    119 &  4 &           &      N &        none & 78 &  157 & 5 & \cite[I]{lu-darmofal04gaussian}  \\

    \end{tabular}
  \end{center}
    {*} A rule with fewer points was known.
\end{table}

\begin{table}
  \caption{ 21 New Cubature Rules for $E_n^r$.}\label{tab:summaryE1}
  \begin{center}
    \begin{tabular}{ccc @{\hspace{1ex}}c@{\hspace{-2pt}} c@{\hspace{2ex}}ccccc}
      \multicolumn{6}{c}{New Rule} &
      \multirow{2}{*}{$N_{MLB}$}&
      \multicolumn{3}{c} {Smallest Previous Rule} \\
      \cmidrule(lr){1-6}  \cmidrule(lr){8-10}
      $n$ & $N$ & $d$ & Shells & Quality & Symmetry & & $N$ & $d$ & Source \\
      \cmidrule(lr){1-6} \cmidrule(lr){7-7}  \cmidrule(lr){8-10}
        2 &     11 &  6 &           &      P &   bilateral & 10 &   12 & 7 & \cite[VI]{stroud-secrest63spherically}  \\
  2 &     17 &  8 &           &      P &   bilateral & 15 &   19 & 9 & \cite{haegemans-piessens77seven}  \\
  3 &     11 &  4 &           &      P &   bilateral & 10 &   13 & 5 & \cite[VII]{stroud-secrest63spherically}  \\
  3 &     23 &  6 &           &      P &        none & 20 &   27 & 7 & \cite[$E_n^r$:7-1]{stroud71approximate}  \\
  4 &     16 &  4 &    1+10+5 &      P &   4-simplex & 15 &   24 & 5 & \cite[$E_n^r$:5-2]{stroud71approximate}  \\
  4 &     16 &  4 &     1+6+9 &      P &   $x_1,x_3$ & 15 &   24 & 5 & \cite[$E_n^r$:5-2]{stroud71approximate}  \\
  4 &     23 &  5 &      1+22 &      P &      vertex & 21 &   24 & 5 & \cite[$E_n^r$:5-2]{stroud71approximate}  \\
  4 &     45 &  6 &           &      P &        none & 35 &   49 & 7 & \cite[$E_n^r$:7-1]{stroud71approximate}  \\
  4 &    103 &  8 &           &      N &        none & 70 &\multicolumn{3}{c}{none} \\
  4 &    154 &  9 &           &      N &        none & 91 &\multicolumn{3}{c}{none} \\
  5 &     22 &  4 &    1+6+15 &      P &   5-simplex & 21 &   42 & 5 & \cite[$E_n^r$:5-2]{stroud71approximate}  \\
  5 &     80 &  6 &           &      P &        none & 56 &   83 & 7 & \cite[$E_n^r$:7-1]{stroud71approximate}  \\
  5 &    230 &  8 &           &      N &        none &126 &\multicolumn{3}{c}{none} \\
  6 &     28 &  4 &      1+27 &      P &      vertex & 28 &   57 & 5 & \cite{meng-luo13few}  \\
  6 &     44 &  5 &     12+32 &      P &     central & 43 &   57 & 5 & \cite{meng-luo13few}  \\
  7 &     46 &  4 &           &      P &        none & 36 &   99 & 5 & \cite[$E_n^r$:5-1]{stroud71approximate}  \\
  7 &    223 &  6 &           &      P &        none &120 &  227 & 7 & \cite[$E_n^r$:7-1]{stroud71approximate}  \\
  8 &     59 &  4 &           &      P &        none & 45 &  129 & 5 & \cite[$E_n^r$:5-1]{stroud71approximate}  \\
  9 &     78 &  4 &           &      P &        none & 55 &  163 & 5 & \cite[$E_n^r$:5-1]{stroud71approximate}  \\
 10 &    107 &  4 &           &      P &        none & 66 &  201 & 5 & \cite[$E_n^r$:5-1]{stroud71approximate}  \\
 11 &    133 &  4 &           &      P &        none & 78 &  243 & 5 & \cite[$E_n^r$:5-1]{stroud71approximate}  \\

    \end{tabular}
  \end{center}
\end{table}

\begin{table}
  \caption{ 21 New Cubature Rules for $S_n$.}\label{tab:summaryS}
  \begin{center}
    \begin{tabular}{ccc @{\hspace{1ex}}c@{\hspace{-4pt}} c@{\hspace{2ex}}ccccc}
     
      \multicolumn{6}{c}{New Rule} &
      \multirow{2}{*}{$N_{MLB}$}&
      \multicolumn{3}{c} {Smallest Previous Rule} \\
      \cmidrule(lr){1-6}  \cmidrule(lr){8-10}
      $n$ & $N$ & $d$ & Shells & Quality & Symmetry & & $N$ & $d$ & Source \\
      \cmidrule(lr){1-6} \cmidrule(lr){7-7}  \cmidrule(lr){8-10}
        2 &    17\clap{*} &  8 &           &     PO &   bilateral & 15 &   16 & 8 & \cite{wissmann-becker86cubature}  \\
  2 &     23 & 10 &           &     PO &   bilateral & 21 &   25 &11 & \cite{haegemans-piessens76eleven}  \\
  3 &     10 &  4 &           &     PO &   $x_2,x_3$ & 10 &   11 & 4 & \cite{becker87finite}  \\
  3 &     22 &  6 &           &     PO &   bilateral & 20 &   27 & 7 & \cite[$S_n$:7-1]{stroud71approximate}  \\
  3 &     42 &  8 &           &     PO &        none & 35 &   45 & 8 & \cite{cools-rabinowitz93cubature}  \\
  4 &16\clap{*}&4 &    1+10+5 &     NO &   4-simplex & 15 &16 & 4 & \cite{cools99cubature,mysovskih92cubature}  \\
  4 &     15 &  4 &     0+6+9 &     PB &   $x_1,x_3$ & 15 &   16 & 4 & \cite{cools99cubature,mysovskih92cubature}  \\
  4 &23\clap{*}&5 &      1+22 &     PI &      vertex & 21 &   22 & 5 & \cite[$S_n$:5-1]{stroud71approximate}  \\
  4 &     43 &  6 &           &     NO &   bilateral & 35 &   49 & 7 & \cite[$S_n$:7-1]{stroud71approximate}  \\
  4 &    105 &  8 &           &     NO &        none & 70 &  193 & 9 & \cite[$S_n$:9-1]{stroud71approximate}  \\
  4 &    147 &  9 &           &     NO &        none & 91 &  193 & 9 & \cite[$S_n$:9-1]{stroud71approximate}  \\
  4 &    208 & 10 &           &     NO &   bilateral &126 &  417 &11 & \cite[$S_n$:11-1]{stroud71approximate}  \\
  5 &     22 &  4 &    1+6+15 &     PI &   5-simplex & 21 &   32 & 5 & \cite[$S_n$:5-1]{stroud71approximate}  \\
  5 &     80 &  6 &           &     NO &        none & 56 &   83 & 7 & \cite[$S_n$:7-1]{stroud71approximate}  \\
  5 &    220 &  8 &           &     NO &        none &126 &  421 & 9 & \cite[$S_n$:9-1]{stroud71approximate}  \\
  6 &     28 &  4 &      1+27 &     PI &      vertex & 28 &   44 & 5 & \cite[$S_n$:5-1]{stroud71approximate}  \\
  6 &    127 &  7 &   1+54+72 &     PB &     central &124 &  137 & 7 & \cite[$S_n$:7-1]{stroud71approximate}  \\
  7 &     38 &  4 &  1+8+8+21 &     NO &    see text & 36 &   57 & 5 & \cite[$S_n$:5-1]{stroud71approximate}  \\
  7 &    183 &  7 &  1+56+126 &     PI &     central &182 &  227 & 7 & \cite[$S_n$:7-1]{stroud71approximate}  \\
  9 &     78 &  4 &           &     NO &        none & 55 &  163 & 5 & \cite[$S_n$:5-2]{stroud71approximate}  \\
 10 &     96 &  4 &           &     NO &        none & 66 &  201 & 5 & \cite[$S_n$:5-2]{stroud71approximate}  \\
 11 &    123 &  4 &           &     NO &        none & 78 &  243 & 5 & \cite[$S_n$:5-2]{stroud71approximate}  \\

    \end{tabular}
  \end{center}
    {*} A rule with fewer or equal points was known.
\end{table}
\subsection{Examples}\label{examples}
To illustrate the formulas, we numerically evaluate an integral used
as an example by Stroud \cite{stroud66fifth}:
\begin{equation}
  \label{eq:ex2}
  J_4=\int_{S_4} \cos(x_1+ \cdots + x_4) \; dx_1 \cdots dx_4 =3.4823322817\;.
\end{equation}
The values calculated using our seven formulas of dimension 4, plus
the 16 point formula of Mysovskih \cite{cools99cubature,mysovskih92cubature},
the 31 point formula of degree 5 by Meng and Luo \cite{meng-luo13few},
and Stroud's formulas of degrees 5 \cite{stroud67fifth} and 7
\cite{stroud67seventh}
are shown
in Table \ref{tab:ex2}, in order by $N$.  
\begin{table}
  \caption{Approximate Values of $J_4$ in \eqref{eq:ex2}.}\label{tab:ex2}
  \centering
  \begin{tabular}{crccrc}
    $n$&$N$ & $d$ & $J_4$ Estimates & Error\hspace{18pt} & Source\\
    \hline
    4 &  15 &  4 & 3.4818127309 &  -0.0005195508 & Table \ref{tab:4_16_4b}\\ 
    4 &  16 &  4 & 3.4511488638 &  -0.0311834178 & Table \ref{tab:4_16_4a}\\
    4 &  16 &  4 & 3.4828928259 &   0.0005605442 & \cite{cools99cubature,mysovskih92cubature}\\
    4 &  22 &  5 & 3.4403244866 &  -0.0420077951 & \cite{stroud67fifth},\cite[$S_n$:5-1]{stroud71approximate}\\
    4 &  23 &  5 & 3.4838622252 &   0.0015299435 & Table \ref{tab:4_23_5}\\
    4 &  31 &  5 & 3.4827186240 &   0.0003863423 & \cite{meng-luo13few}\\
    4 &  43 &  6 & 3.4823547183 &   0.0000224367 & Table \ref{tab:summaryS}\\
    4 &  49 &  7 & 3.4823164472 &  -0.0000158345 & \cite{stroud67seventh},\cite[$S_n$:7-1]{stroud71approximate}\\
    4 & 105 &  8 & 3.4823287423 &  -0.0000035394 & Table \ref{tab:summaryS}\\
    4 & 147 &  9 & 3.4823311982 &  -0.0000010835 & Table \ref{tab:summaryS}\\
    4 & 208 & 10 & 3.4823322804 &  -0.0000000012 & Table \ref{tab:summaryS}\\
  \end{tabular}
\end{table}

\subsection{Supplemental material}
The supplemental material includes plain-text listings of the new
rules, in double and quad precision, with 15 and 32 decimal digits,
respectively.  Some known rules are included for comparison, with
sources indicated in the double precision listings.  
The quad precision listings are of two sorts, both
generated by programs in Maxima.  Where
closed form expressions were found for the parameters of a rule, those
expressions were evaluated with 64 digit precision and printed
with 32 digit precision.  Otherwise, a simple root-finder using
Newton's method with Moore-Penrose pseudoinverses
was used to
refine the double precision rule with an excessive $32d+10$ digits of
precision, with the goal that the printed values of both node
coordinates and weights would be correct to 32 digits.
In either case, the constraint equations to the stated
degree were evaluated and the maximum error printed.  The error is
zero where the parameters were expressed in closed form and Maxima was
able to simplify the resulting equations.  Otherwise the error is the result of an
extended precision calculation.

Also included are
several of the MATLAB/Octave and Maxima programs used to find these
rules and to refine them to high precision.

%

\section*{Acknowledgements}
I thank one anonymous reviewer for finding closed forms for several
coefficients in the rule of degree 4 for $S_3$, and another anonymous
reviewer for providing the weights and points for Mysovskih's 16 point
formula for $S_4$.


%
%


\ifels
\bibliographystyle{elsarticle-num}
\fi
 
  \bibliography{sigma}
\end{document}